\documentclass{amsart}
\usepackage{amssymb}

\copyrightinfo{2001}{American Mathematical Society}

\newtheorem{theorem}{Theorem}

\theoremstyle{definition}

\theoremstyle{remark}

\begin{document}

\title{Generators and relations for Schur algebras}

\author{Stephen Doty}
\address{Loyola University Chicago, Chicago, Illinois 60626 U.S.A.}
\email{doty@math.luc.edu}

\author{Anthony Giaquinto}
\address{Loyola University Chicago, Chicago, Illinois 60626 U.S.A.}
\email{tonyg@math.luc.edu}



\date{April 2, 2001.}





\newcommand{\N}{{\mathbb N}}
\newcommand{\Z}{{\mathbb Z}}
\newcommand{\Q}{{\mathbb Q}}
\newcommand{\R}{{\mathbb R}}
\newcommand{\C}{{\mathbb C}}
\newcommand{\g}{{\mathfrak g}}
\newcommand{\n}{{\mathfrak n}}
\newcommand{\dist}{\operatorname{Dist}}
\newcommand{\per}{\operatorname{per}}
\newcommand{\cov}{\operatorname{cov}}
\newcommand{\non}{\operatorname{non}}
\newcommand{\cf}{\operatorname{cf}}
\newcommand{\add}{\operatorname{add}}
\newcommand{\End}{\operatorname{End}}
\newcommand{\Ext}{\operatorname{Ext}}
\newcommand{\Hom}{\operatorname{Hom}}
\newcommand{\Tor}{\operatorname{Tor}}
\newcommand{\ch}{\operatorname{ch}}
\newcommand{\ind}{\operatorname{ind}}
\newcommand{\coind}{\operatorname{Coind}}
\newcommand{\res}{\operatorname{res}}
\newcommand{\soc}{\operatorname{soc}}
\newcommand{\rad}{\operatorname{rad}}
\newcommand{\Aut}{\operatorname{Aut}}
\newcommand{\Dist}{\operatorname{Dist}}
\newcommand{\Lie}{\operatorname{Lie}}
\renewcommand{\ker}{\operatorname{Ker}}
\newcommand{\im}{\operatorname{im}}

\newcommand{\GL}{{\sf GL}}
\newcommand{\SL}{{\sf SL}}
\newcommand{\gl}{\mathfrak{gl}}
\renewcommand{\sl}{\mathfrak{sl}}
\newcommand{\B}{\mathcal{B}}
\newcommand{\divided}[2]{#1^{(#2)}}
\newcommand{\ct}{\chi}
\newcommand{\ctl}{\ct_L}
\newcommand{\ctr}{\ct_R}
\newcommand{\sqbinom}[2]{\begin{bmatrix}#1\\#2\end{bmatrix}}
\newcommand{\U}{\mathbf{U}}
\renewcommand{\S}{\mathbf{S}}
\newcommand{\A}{\mathcal{A}}


\begin{abstract}
We obtain a presentation of Schur algebras (and $q$-Schur algebras) by
generators and relations, one which is compatible with the usual
presentation of the enveloping algebra (quantized enveloping algebra)
corresponding to the Lie algebra $\gl_n$ of $n\times n$ matrices.  We
also find several new bases of Schur algebras.
\end{abstract}

\maketitle

\parskip=2pt

\section{Introduction}\label{sec:intro}
The classical Schur algebra $S(n,d)$ (over $\Q$) may be defined as the
algebra $\End_{\Sigma_d}(V^{\otimes d})$ of linear endomorphisms on
the $d$th tensor power of an $n$-dimensional $\Q$-vector space $V$
commuting with the action of the symmetric group $\Sigma_d$, acting by
permutation of the tensor places (see \cite{Green}).  Schur algebras
determine the polynomial representation theory of general linear
groups, and they form an important class of quasi-hereditary algebras.

We identify $V$ with $\Q^n$. Then $V_\Z=\Z^n$ is a lattice in
$V$.  One can define a $\Z$-order $S_\Z(n,d)$ (the integral Schur
algebra) in $S(n,d)$ by setting $S_\Z(n,d) =
\End_{\Sigma_d}(V_\Z^{\otimes d})$.  For any field $K$, one then
obtains the Schur algebra $S_K(n,d)$ over $K$ by setting $S_K(n,d)=
S_\Z(n,d)\otimes_\Z K$.  Moreover, $S_\Q(n,d) \cong S(n,d)$.

In the quantum case one can replace $\Q$ by $\Q(v)$ ($v$ an
indeterminate), $V$ by an $n$-dimensional $\Q(v)$-vector space, and
$\Sigma_d$ by the corresponding Hecke algebra $\mathbf{H}(\Sigma_d)$.
Then the resulting commuting algebra, $\S(n,d)$, is known as the
$q$-Schur algebra, or quantum Schur algebra. It appeared first in work
of Dipper and James, and, independently, Jimbo. Beilinson, Lusztig,
and MacPherson \cite{BLM} have given a geometric construction of
$\S(n,d)$ in terms of orbits of flags in vector spaces. (See also
\cite{Du}.)

In this situation $\Z$ is replaced by the ring $\A=\Z[v,v^{-1}]$, and
there is a corresponding ``integral'' form $\S_\A(n,d)$ in
$\S(n,d)$. The above quantized objects specialize to their classical
versions when $v=1$.  We have more detailed information about the
Schur algebras in rank $1$; see \cite{DG} and \cite{DG:quantum}.
Proofs of the main results will appear in \cite{PSA}.

\section{Serre's presentation of $U$}\label{sec:presentU}

Let $\Phi=\{\varepsilon_i - \varepsilon_j \mid 1\le i \ne j \le n\}$
be the root system of type $A_{n-1}$, with simple roots $\Delta=\{
\varepsilon_i - \varepsilon_{i+1} \mid i=1, \dots, n-1 \}$, where
$\{\varepsilon_1, \dots, \varepsilon_n\}$ denotes the standard basis
of $\Z^n$.  The corresponding set $\Phi^+$ of positive roots is the
set $ \Phi^+ = \{ \varepsilon_i - \varepsilon_j \mid 1\le i < j \le n
\}$.

The abelian group $\Z^n$ has a bilinear form $(\ ,\ )$ given by
$(\varepsilon_i,\varepsilon_j)=\delta_{ij}$ (Kronecker delta). We
write $\alpha_j:= \varepsilon_j - \varepsilon_{j+1}$.

The enveloping algebra $U=U(\gl_n)$ is the associative algebra (with
$1$) on generators $e_i,f_i$ ($i=1,\dots,n-1$) and $H_i$
($i=1,\dots,n$) with relations
\begin{align}
H_i H_j &= H_j H_i \tag{R1}\label{R1} \\
e_i f_j - f_j e_i &= \delta_{ij}(H_j - H_{j+1})  \tag{R2}\label{R2}\\
{}H_i e_j - e_j H_i = (\varepsilon_i, \alpha_j) e_j&, \quad
H_i f_j - f_j H_i = -(\varepsilon_i, \alpha_j) f_j \tag{R3}\label{R3}
\end{align}
\begin{equation}\label{R4}
\begin{aligned}
{}&e_i^2 e_j - 2 e_i e_j e_i + e_j e_i^2 = 0 \quad (|i-j|=1)\\
&e_i e_j - e_j e_i = 0 \quad (\text{otherwise})
\end{aligned} \tag{R4}
\end{equation}
\begin{equation}\label{R5}
\begin{aligned}
{}&f_i^2 f_j - 2 f_i f_j f_i + f_j f_i^2 = 0 \quad (|i-j|=1)\\
&f_i f_j - f_j f_i = 0 \quad (\text{otherwise}).
\end{aligned} \tag{R5}
\end{equation}

For $\alpha\in \Phi$, let $x_{\alpha}$ denote the corresponding
root vector of $\gl_n$ viewed as an element of $U$. We have in
particular $e_i=x_{\alpha_i}$ and $f_i=x_{-\alpha_i}.$

\newcommand{\ep}{\varepsilon}
\newcommand{\al}{\alpha}
\newcommand{\be}{\beta}
\newcommand{\Ea}{E_{\alpha}}
\newcommand{\Eb}{E_{\beta}}
\newcommand{\Eg}{E_{\gamma}}
\newcommand{\Ep}{E_{\rho}}
\newcommand{\Eij}{E_{ij}}
\newcommand{\Ekl}{E_{kl}}
\newcommand{\vbinom}[2]{\begin{bmatrix}#1\\#2\end{bmatrix}}
\newcommand{\vbinomm}[3]{\begin{bmatrix}#1; #2\\#3\end{bmatrix}}
\newcommand{\1}{1_{\lambda}}
\newcommand{\Lnd}{\Lambda (n,d)}
\newcommand{\Ua}{\mathbf{U}_{\mathcal{A}}}

\section{The quantized enveloping algebra}\label{sec:que}

The Drinfeld-Jimbo quantized enveloping algebra
$\U=\U_v(\gl_n)$, by definition, is the $\Q(v)$-algebra with
generators $E_i$, $F_i$ ($1\leq i\leq n-1$), $K_i^{\pm
1}$ ($1\leq i\leq n$) and relations
\begin{align}
& K_i K_j =K_j K_i, \qquad K_i K_i^{-1} = K_i^{-1}K_i =1
  \tag{Q1}\label{Q1}\\
& E_iF_j - F_jE_i =
\delta_{ij}\frac{K_iK_{i+1}^{-1}-K_i^{-1}K_{i+1}}{v-v^{-1}}\tag{Q2}\label{Q2}\\
& K_i E_j = v^{(\ep_i, \al_j)}E_jK_i, \qquad
K_i F_j = v^{-(\ep_i, \al_j)}F_jK_i \tag{Q3}\label{Q3}
\end{align}
\begin{equation}\label{Q4}
\begin{aligned}
{}& E_i^2E_j - (v+v^{-1})E_iE_jE_i +E_jE_i^2 =0 \quad(|i-j|=1) \\
{}& E_iE_j - E_i E_j = 0 \quad(\text{otherwise})
\end{aligned}\tag{Q4}
\end{equation}
\begin{equation}
\begin{aligned}\label{Q5}
{}& F_i^2F_j - (v+v^{-1})F_iF_jF_i +F_jF_i^2 =0 \quad(|i-j|=1) \\
{}& F_iF_j - F_jF_i = 0 \quad(\text{otherwise}).
\end{aligned}\tag{Q5}
\end{equation}

For $\alpha\in \Phi^+$, let $E_{\alpha}$ and $F_{\alpha}$ be the
positive and negative quantum root vectors of $\U$ as defined by
Jimbo \cite{Jimbo}. We have in particular $E_i=E_{\alpha_i}$ and
$F_i=F_{\alpha_i}.$

\section{Main results: classical case} \label{sec:main}
We now give a precise statement of our main results in the classical
case.  The first result describes a presentation by generators and
relations of the Schur algebra over the rational field $\Q$. This
presentation is compatible with the usual presentation (see section
\ref{sec:presentU}) of $U=U(\gl_n)$.

\begin{theorem} \label{thm:present:S}
Over $\Q$, the Schur algebra $S(n,d)$ is isomorphic to the
associative algebra (with 1) on the same generators as for $U$
subject to the same relations \eqref{R1}--\eqref{R5} as for $U$,
together with the additional relations:
\begin{align}
& H_1 + H_2 + \cdots + H_n = d \tag{R6}\label{R6} \\
& H_k(H_k-1)\cdots(H_k-d) = 0, \qquad(k = 1, \cdots, n).
  \tag{R7}\label{R7}
\end{align}
\end{theorem}

The next result gives a basis for the Schur algebra which is the analogue
of the Poincare-Birkhoff-Witt (PBW) basis of $U$.

\begin{theorem} \label{thm:truncPBW}
The algebra $S(n,d)$ has a ``truncated PBW'' basis (over $\Q$) which
can be described as follows. Fix an integer $k_0$ with $1\le k_0 \le
n$ and set
$$
G=\{ x_\alpha \mid \alpha\in \Phi \} \cup
\{H_k \mid k \in \{1,\dots,n\}-\{k_0\}\, \}
$$
and fix some ordering for this set.  Then the set of all monomials in
$G$ (with specified order) of total degree at most $d$ is a basis for
$S(n,d)$.
\end{theorem}

Our next result constructs the integral Schur algebra $S_\Z(n,d)$
in terms of the generators given above.  We need some more
notation. For $B$ in $\N^n$, we write
$$
H_B = \prod_{k=1}^n \dbinom{H_k}{b_k}.
$$
Let $\Lambda(n,d)$ be the subset of $\N^n$ consisting of those
$\lambda \in \N^n$ satisfying $|\lambda|=d$; this is the set of
$n$-part compositions of $d$. Given $\lambda \in \Lambda(n,d)$ we
set $1_\lambda:=H_{\lambda}.$ One can show that the collection
$\{1_\lambda\}$ as $\lambda$ varies over $\Lambda(n,d)$ forms a
set of pairwise orthogonal idempotents in $S_\Z(n,d)$ which sum
to the identity element.

For $m\in \N$ and $\alpha \in \Phi$, set $\divided{x_\alpha}{m}:=
x_\alpha^m/(m!)$. Any product in $U$ of elements of the form
$$
\divided{x_\alpha}{m}, \quad \dbinom{H_k}{m}\qquad
(m\in \N, \alpha\in\Phi, k\in \{1,\dots, n\}),
$$
taken in any order, will be called a {\em Kostant monomial}.  Note
that the set of Kostant monomials is multiplicatively closed, and
spans $U_\Z$.  We define a function $\ct$ (content function) on
Kostant monomials by setting
$$
\ct(\divided{x_\alpha}{m}):= m\,\varepsilon_{\max(i,j)}, \quad
\ct( \binom{H_k}{m} ) := 0
$$
where $\alpha = \varepsilon_i - \varepsilon_j$ ($i\ne j$), and by
declaring that $\ct(XY) = \ct(X)+\ct(Y)$ whenever $X,Y$ are Kostant
monomials.

We write any $A\in \N^{\Phi^+}$ in ``multi-index'' form
$A=(a_\alpha)_{\alpha\in \Phi^+}$ and set $|A|:= \sum_{\alpha\in
\Phi^+} a_\alpha$.  For $A, C \in \N^{\Phi^+}$ we write
$$
e_A = \prod_{\alpha\in \Phi^+} \divided{x_\alpha}{a_\alpha}, \quad
f_C = \prod_{\alpha\in \Phi^-} \divided{x_{\alpha}}{c_\alpha}
$$
where the products in $e_A$ and $f_C$ are taken relative to any
fixed orders on $\Phi^+$ and $\Phi^-$.

\begin{theorem} \label{thm:idemp:basis}
The integral Schur algebra $S_\Z(n,d)$ is the subring of $S(n,d)$
generated by all divided powers $ \divided{e_j}{m},\,
\divided{f_j}{m} (j\in\{1,\ldots , n-1\}, m\in \N)$. Moreover,
this algebra has a $\Z$-basis consisting of all
\begin{equation}\label{B1}
e_A 1_\lambda f_C  \qquad (\ct(e_A f_C) \preceq \lambda)\tag{B1}
\end{equation}
and another such basis consisting of all
\begin{equation}\label{B2}
f_A 1_\lambda e_C  \qquad (\ct(f_A e_C) \preceq \lambda)\tag{B2}
\end{equation}
where $A,C \in \N^{\Phi^+}$, $\lambda\in \Lambda(n,d)$, and where
$\preceq$ denotes the {\em componentwise} partial ordering on
$\N^n$.
\end{theorem}

Finally, we have another presentation of the Schur algebra by
generators and relations. This presentation has the advantage that it
possesses a quantization of the same form, in which we can set $v=1$
to recover the classical version.

\begin{theorem} \label{thm:idemp:present}
For each $i\in \{1,\dots,n-1\}$ write $\alpha_i := \varepsilon_i
-\varepsilon_{i+1}$.
The algebra $S(n,d)$ is the associative
algebra (with $1$) given by generators $1_\lambda$ ($\lambda\in
\Lambda(n,d)$), $e_i$, $f_i$ ($i\in \{1,\dots,n-1\}$) subject to
the relations
\begin{equation}\label{S1}
\begin{gathered}
1_\lambda 1_\mu = \delta_{\lambda,\mu} 1_\lambda
  \quad(\lambda,\mu\in \Lambda(n,d))\\
\sum_{\lambda\in \Lambda(n,d)} 1_\lambda = 1
\end{gathered}\tag{S1}
\end{equation}
\begin{equation}\label{S2}
\begin{gathered}
e_i 1_\lambda =
\begin{cases}
1_{\lambda+\alpha_i} e_i &
   \text{if $\lambda+\alpha_i \in \Lambda(n,d)$}\\
0 & \text{otherwise}
\end{cases}\\
f_i 1_\lambda =
\begin{cases}
1_{\lambda-\alpha_i} f_i &
   \text{if $\lambda-\alpha_i \in \Lambda(n,d)$}\\
0 & \text{otherwise}
\end{cases}\\
1_\lambda e_i =
\begin{cases}
e_i 1_{\lambda-\alpha_i} &
   \text{if $\lambda-\alpha_i \in \Lambda(n,d)$}\\
0 & \text{otherwise}
\end{cases}\\
1_\lambda f_i =
\begin{cases}
f_i 1_{\lambda+\alpha_i} &
   \text{if $\lambda+\alpha_i \in \Lambda(n,d)$}\\
0 & \text{otherwise}
\end{cases}
\end{gathered}\tag{S2}
\end{equation}
\begin{equation}\label{S3}
e_i f_j - f_j e_i = \delta_{ij}
\sum_{\lambda\in \Lambda(n,d)} (\lambda_j-\lambda_{j+1}) 1_\lambda \tag{S3}
\end{equation}
along with the Serre relations \eqref{R4}, \eqref{R5}, for
$i,j \in \{1, \dots,n-1\}$.
\end{theorem}

\section{Main results: quantized case} \label{qsec:main}
Our main results in the quantized case are similar in form to those in
the classical case.  The first result describes a presentation by
generators and relations of the quantized Schur algebra over the
rational function field $\Q(v)$. This presentation is compatible with
the usual presentation (see section \ref{sec:que}) of $\U=\U_v(\gl_n)$.

\setcounter{theorem}{0} \renewcommand{\thetheorem}{$\arabic{theorem}^\prime$}

\begin{theorem} \label{qthm:present:S}
Over $\Q(v)$, the $q$-Schur algebra $\S(n,d)$ is isomorphic with the
associative algebra (with 1) on the same generators as for $\U$
subject to the same relations \eqref{Q1}--\eqref{Q5} as
for $\U$, together with the additional relations:
\begin{align}
& K_1  K_2 \cdots  K_n = v^d \tag{Q6}\label{Q6} \\
& (K_j-1)(K_j-v)(K_j-v^2)\cdots(K_j-v^d) = 0, \qquad(j=1, \cdots, n).
  \tag{Q7}\label{Q7}
\end{align}
\end{theorem}

The next result gives a basis for the $q$-Schur algebra which is
the analogue of the Poincare-Birkhoff-Witt (PBW) type basis of
$\U$, given in Lusztig \cite[Prop.\ 1.13]{Lusztig:book}.

\begin{theorem} \label{qthm:truncPBW}
The algebra $\S(n,d)$ has a ``truncated PBW type'' basis which
can be described as follows. Fix an integer $k_0$ with $1\le k_0
\le n$ and set
$$
G'=\{ E_\alpha, F_{\alpha} \mid \alpha\in \Phi^+ \} \cup \{K_k
\mid k \in \{1,\dots,n\}-\{k_0\}\, \}
$$
and fix some ordering for this set.  Then the set of all
monomials in $G'$ (with specified order) of total degree at most
$d$ is a basis for $\S(n,d)$.
\end{theorem}

Note that setting $v=1$ in the basis of Theorem
\ref{qthm:truncPBW} does not yield the basis of Theorem
\ref{thm:truncPBW} since $K_i$ acts as the identity when $v=1$.

Our next result constructs the integral $q$-Schur algebra
$\S_\A(n,d)$ in terms of the generators given above. For $B$ in
$\N^n$, we write
$$
K_B = \prod_{j=1}^n \sqbinom{K_j}{b_j}.
$$
where for indeterminates $X,X^{-1}$ satisfying $X X^{-1} = X^{-1} X = 1$ and
any $t\in \N$ we formally set
$$
\sqbinom{X}{t}:=\prod_{s=1}^{t}
\frac{Xv^{-s+1}-X^{-1}v^{s-1}}{v^s-v^{-s}},
$$
an expression that makes
sense if $X$ is replaced by any invertible element of a
$\Q(v)$-algebra.

Given $\lambda \in \Lambda(n,d)$ we set (when we are in the
quantum case) $1_\lambda:= K_\lambda$. It turns out that, just as
in the classical case, the collection $\{1_\lambda\}$ as
$\lambda$ varies over $\Lambda(n,d)$ forms a set of pairwise
orthogonal idempotents in $\S_\A(n,d)$ which sum to the identity
element.

Let $[m]$ denote the quantum integer $[m]:=
(v^m-v^{-m})/(v-v^{-1})$ and set $[m]! := [m][m-1] \cdots [1]$.
Then the $q$-analogues of the divided powers of root vectors are
defined to be $ \divided{E_\alpha}{m}:= E_\alpha/[m]!$ and
$\divided{F_\alpha}{m}:= F_\alpha/[m]!.$ Any product in $\U$ of
elements of the form
$$
\divided{E_\alpha}{m}, \quad \divided{F_\alpha}{m}, \quad\sqbinom{K_j}{m}
\qquad (m\in \N, \alpha\in\Phi, j\in \{1,\dots, n\}),
$$
taken in any order, will be called a {\em Kostant monomial}.  As before,
the set of Kostant monomials is multiplicatively closed, and
spans $\U_\A$.  The definition of content $\ct$ of
a Kostant monomial is obtained similarly, by setting
$$
\ct(\divided{E_\alpha}{m}) = \ct(\divided{F_\alpha}{m})
  := m\,\varepsilon_{\max(i,j)}, \quad
\ct( \sqbinom{K_l}{m} ) := 0
$$
where $\alpha = \varepsilon_i - \varepsilon_j \in \Phi^+$, and by
declaring that $\ct(XY) = \ct(X)+\ct(Y)$ whenever $X,Y$ are
Kostant monomials. For $A, C \in \N^{\Phi^+}$ we write
$$
E_A = \prod_{\alpha\in \Phi^+} \divided{E_\alpha}{a_\alpha}, \quad
F_C = \prod_{\alpha\in \Phi^+} \divided{F_{\alpha}}{c_\alpha}
$$
where the products in $E_A$ and $F_C$ are taken relative to any
two specified orderings on $\Phi^+$.

\begin{theorem} \label{qthm:idemp:basis}
The integral $q$-Schur algebra $\S_\A(n,d)$ is the subring of
$\S(n,d)$ generated by all quantum divided powers $
\divided{E_j}{m},\, \divided{F_j}{m} \,(j\in \{ 1,\ldots , n-1\},
m\in \N),$ along with the elements $\sqbinom{K_j}{m}$ ($j \in
\{1,\dots,n\}$, $m\in \N$). Moreover, this algebra has a basis
over $\A$ consisting of all
\begin{equation}\label{B1'}
E_A 1_\lambda F_C  \qquad (\ct(e_A f_C) \preceq \lambda)\tag{B$1^\prime$}
\end{equation}
and another such basis consisting of all
\begin{equation}\label{B2'}
F_A 1_\lambda E_C  \qquad (\ct(f_A e_C) \preceq \lambda)\tag{B$2^\prime$}
\end{equation}
where $A,C \in \N^{\Phi^+}$, $\lambda\in \Lambda(n,d)$, and where
$\preceq$ denotes the {\em componentwise} partial ordering on
$\N^n$.
\end{theorem}

Unlike the truncated PBW basis, the bases of Theorem
\ref{qthm:idemp:basis} do specialize when $v=1$ to their classical
analogues given in Theorem \ref{thm:idemp:basis}.

Finally, we have another presentation of the $q$-Schur algebra by
generators and relations. This presentation has the advantage that by
setting $v=1$, we recover the classical version given in Theorem
\ref{thm:idemp:present}.  The relations of the presentation are
similar to relations that hold for Lusztig's modified form $\dot{\U}$
of $\U$.  (See \cite[Chap.\ 23]{Lusztig:book}.)

\begin{theorem} \label{qthm:idemp:present}
For each $i\in \{1,\dots,n-1\}$ write $\alpha_i := \varepsilon_i
-\varepsilon_{i+1}$.
The algebra $\S(n,d)$ is the associative
algebra (with $1$) given by generators $1_\lambda$ ($\lambda\in
\Lambda(n,d)$), $E_i$, $F_i$ ($i\in \{1,\dots,n-1\}$) subject to
the relations
\begin{equation}\label{S1'}
\begin{gathered}
1_\lambda 1_\mu = \delta_{\lambda,\mu} 1_\lambda
  \quad(\lambda,\mu\in \Lambda(n,d))\\
\sum_{\lambda\in \Lambda(n,d)} 1_\lambda = 1
\end{gathered} \tag{S$1^\prime$}
\end{equation}
\begin{equation}\label{S2'}
\begin{gathered}
E_i 1_\lambda =
\begin{cases}
1_{\lambda+\alpha_i} E_i &
   \text{if $\lambda+\alpha_i \in \Lambda(n,d)$}\\
0 & \text{otherwise}
\end{cases}\\
F_i 1_\lambda =
\begin{cases}
1_{\lambda-\alpha_i} F_i &
   \text{if $\lambda-\alpha_i \in \Lambda(n,d)$}\\
0 & \text{otherwise}
\end{cases}\\
1_\lambda E_i =
\begin{cases}
E_i 1_{\lambda-\alpha_i} &
   \text{if $\lambda-\alpha_i \in \Lambda(n,d)$}\\
0 & \text{otherwise}
\end{cases}\\
1_\lambda F_i =
\begin{cases}
F_i 1_{\lambda+\alpha_i} &
   \text{if $\lambda+\alpha_i \in \Lambda(n,d)$}\\
0 & \text{otherwise}
\end{cases}
\end{gathered}\tag{S$2^\prime$}
\end{equation}
\begin{equation} \label{S3'}
E_i F_j - F_j E_i = \delta_{ij}
\sum_{\lambda\in \Lambda(n,d)} [\lambda_j-\lambda_{j+1}] 1_\lambda
  \tag{S$3^\prime$}
\end{equation}
along with the $q$-Serre relations \eqref{Q4}, \eqref{Q5}, for
$i,j \in \{1, \dots,n-1\}$.
\end{theorem}

\section{Other results}\label{sec:other}

\subsection{Triangular decomposition.}
One can define the plus, minus, and zero parts of Schur algebras in
terms of the generators, as follows. (These subalgebras have been
studied before.) The plus part $S^+(n,d)$ (resp., the minus part
$S^-(n,d)$) is the subalgebra of $S(n,d)$ generated by all $x_\alpha$,
$\alpha\in \Phi^+$ (resp., $\alpha\in \Phi^-$). The zero part
$S^0(n,d)$ is the subalgebra generated by all $H_j$, $j=1,\dots, n$.
We also have the Borel Schur algebras $S^{\geqslant0}(n,d)$ (resp.,
$S^{\leqslant0}(n,d)$), the subalgebra generated by $S^+(n,d)$ (resp.,
$S^-(n,d)$) together with $S^0(n,d)$.

The appellations $S_\Z^+(n,d)$, $S_\Z^-(n,d)$, $S_\Z^0(n,d)$,
$S_\Z^{\geqslant0}(n,d)$, $S_\Z^{\leqslant0}(n,d)$ will denote the
intersection of the appropriate algebra from above with the integral
form $S_\Z(n,d)$.

The algebra $S=S(n,d)$ has a triangular decomposition $S=S^+S^0S^-$.
We show that in this decomposition one can permute the three factors
in any order. Moreover, the same result holds over $\Z$.

The zero part $S_\Z^0(n,d)$ is the algebra generated by all
$\dbinom{H_j}{m}$ for $j=1,\dots,n$ and $m \in \N$. We give in
\cite{PSA} a presentation of $S^0(n,d)$ by generators and relations.
In particular, we prove that $H_B = 0$ whenever $|B| > d$ ($B\in
\N^n$) and that the set of (pairwise orthogonal) idempotents
$1_\lambda$, $\lambda\in \Lambda(n,d)$, is a $\Z$-basis of
$S_\Z^0(n,d)$.

$S_\Z^+(n,d)$ (resp., $S_\Z^-(n,d)$) is the algebra generated by all
divided powers $\divided{x_\alpha}{m}$ for $\alpha\in \Phi^+$ (resp.,
$\alpha\in \Phi^-$) and $m\in \N$.  It is an easy consequence of the
commutation formulas \eqref{S2} that each generator $x_\alpha$
is nilpotent of index $d+1$; see \cite{PSA} for details.  Moreover,
from our main results we see easily that the set of all $e_A$ (resp.,
$f_A$) such that $|A|\le d$ is a $\Z$-basis for the algebra
$S_\Z^+(n,d)$ (resp., $S_\Z^-(n,d)$).

It also follows immediately from our results that the set of all $e_A
1_\lambda$ (resp., $1_\lambda f_A$) satisfying $\ct(e_A)\preceq
\lambda$ is an integral basis for the Borel Schur algebra
$S_\Z^{\geqslant0}(n,d)$ (resp., $S_\Z^{\leqslant0}(n,d)$).

Similar statements to the above hold in the quantum case. In
particular, as a corollary of the commutation formulas
\eqref{S2'} one can give a simple proof of \cite[Prop.\ 2.3]{RGreen}.

\subsection{Explicit reduction formulas.}
Fix a positive root $\alpha$ and write $\alpha = \varepsilon_i -
\varepsilon_j$ for $i<j$.

\newcommand{\Ha}{{H_i}}
\newcommand{\Hb}{{H_j}}

Then from \cite{DG} we have the following reduction formulas in
$S_\Z(n,d)$, for any $a,b,c\in \N$:
\begin{align}
\divided{f_\alpha}{a} \binom{\Hb}{b} \divided{e_\alpha}{c}
&= \sum_{k=s}^{\min(a,c)}(-1)^{k-s}\binom{k-1}{s-1}\binom{b+k}{k}
\divided{f_\alpha}{a-k} \binom{\Hb}{b+k}
\divided{e_\alpha}{c-k}\\
\divided{e_\alpha}{a}\binom{\Ha}{b}\divided{f_\alpha}{c}
&= \sum_{k=s}^{\min(a,c)}
(-1)^{k-s}\binom{k-1}{s-1}\binom{b+k}{k}\divided{e_\alpha}{a-k}
\binom{\Ha}{b+k}\divided{f_\alpha}{c-k}
\end{align}
where $s=a+b+c-d$ and $s\ge 1$.

We do not have a $q$-analogue of these formulas.

\newcommand{\Ka}{{K_1}}
\newcommand{\Kb}{{K_2}}
\newcommand{\Kaa}{K_{1}^{-1}}
\newcommand{\Kbb}{K_{2}^{-1}}
\newcommand{\Ki}{{K_i}}
\newcommand{\KKa}{{\overline{K_1}}}
\newcommand{\KKb}{{\overline{K_2}}}
\newcommand{\vv}{v^{-1}}
\newcommand{\K}{\begin{bmatrix}K_i\\b_1\end{bmatrix}\begin{bmatrix}K_j\\b_2\end{bmatrix}}
\newcommand{\KK}[2]{\begin{bmatrix}K_i\\#1 \end{bmatrix}\begin{bmatrix}K_j\\#2\end{bmatrix}}
\newcommand{\ea}{E_\alpha^{(a)}}
\newcommand{\fc}{F_\alpha^{(c)}}
\newcommand{\fa}{F_\alpha^{(a)}}
\newcommand{\ec}{E_\alpha^{(c)}}

The results of \cite{DG:quantum} give another type of reduction
formula for $\S_\A(n,d)$. If $b_1, b_2\in \N$ satisfy
$b_1+b_2=d$, set $\lambda:=b_1 \varepsilon_i + b_2 \varepsilon_j
\in \Lambda(n,d)$. Then for all $s\ge1$ we have:
\begin{align}\ea 1_{\lambda}\fc&=
\sum_{k=s}^{\min(a,c)}(-1)^{k-s}\sqbinom{k-1}{s-1}\sqbinom{b_1+k}{k}
\divided{E_\alpha}{a-k}1_{\lambda+k\alpha}\divided{F_\alpha}{c-k}\label{qfirst}\\
\fa 1_{\lambda} \ec&=
\sum_{k=s}^{\min(a,c)}(-1)^{k-s}\sqbinom{k-1}{s-1}\sqbinom{b_2+k}{k}
\divided{F_\alpha}{a-k}1_{\lambda-k\alpha}\divided{E_\alpha}{c-k}\label{qsecond}
\end{align}
where $s=a+b_1+c-d$ in \eqref{qfirst} and $s=a+b_2+c-d$ in
\eqref{qsecond}.

The classical analogues of formulas \eqref{qfirst} and
\eqref{qsecond} hold in $S(n,d)$.

\subsection{Another presentation (for $n=2$).}
We have the following result from \cite{DG}, which presents $S(2,d)$
as a quotient of $U(\sl_2)$.

\renewcommand{\thetheorem}{\arabic{theorem}}

\begin{theorem}\label{thm:sl2}
Over $\Q$, the Schur algebra $S(2,d)$ is isomorphic
with the associative algebra (with 1) generated by $e$, $f$, $h$
subject to the relations:
\begin{align*}
& he-eh=2e; \qquad ef-fe=h; \qquad hf-fh=-2f\\
& (h+d)(h+d-2)\cdots(h-d+2)(h-d)=0.
\end{align*}
Moreover, this algebra has a ``truncated PBW'' basis
over $\Q$ consisting of all
$f^a h^b e^c$ such that $a+b+c \le d$.
\end{theorem}

Note that if we eliminate the last relation we have the usual
presentation of $U(\sl_2)$ over $\Q$.  The last relation is the
minimal polynomial of $h$ in the representation on tensor space.  The
problem of presenting $S(n,d)$ as a quotient of $U(\sl_n)$ seems to be
more difficult for $n>2$.

We also have from \cite{DG:quantum} the following $q$-version of the
above, which presents the $q$-Schur algebra $\S(2,d)$ as a quotient of
the quantized enveloping algebra $\U(\sl_2)$.

\setcounter{theorem}{4} \renewcommand{\thetheorem}{$\arabic{theorem}^\prime$}

\begin{theorem}
Over $\Q(v)$, the quantum Schur algebra $\S(2,d)$ is isomorphic to the
algebra generated by $E$, $F$, $K^{\pm 1}$ subject to the
relations:
\begin{align*}
& KK^{-1}=K^{-1}K=1\\
&K E K^{- 1} =v^{2}E \qquad
K F K^{- 1} =v^{- 2}F \\
& EF-FE=\frac{K - K^{-1}}{v-v^{-1}} \\
& (K-v^d)(K-v^{d-2})\cdots (K-v^{-d+2})(K-v^{-d})=0.
\end{align*}
\end{theorem}

\subsection{Hecke algebras.}
Suppose that $n\ge d$.  Let $\omega = (1^d)$.  Then the subalgebra
$1_\omega \S(n,d) 1_\omega$ is isomorphic with the Hecke algebra
$\mathbf{H}(\Sigma_d)$.  The nonzero elements of the basis \eqref{B1'}
of the form $1_\omega E_A F_C 1_\omega$ form a basis of the Hecke
algebra; similarly for elements of the basis \eqref{B2'} of the form
$1_\omega F_A E_C 1_\omega$.

Moreover, taking $d=n$, we can see that
$\mathbf{H}=\mathbf{H}(\Sigma_n)$ is generated by the elements
$1_\omega E_i F_i 1_\omega$ ($1 \le i \le n-1$).  Alternatively,
$\mathbf{H}$ is generated by the $1_\omega F_i E_i 1_\omega$ ($1 \le i
\le n-1$).  One can easily describe the relations on these generators,
thus obtaining a presentation of $\mathbf{H}$ which is closely related
to one in \cite{Wenzl}.


\bibliographystyle{amsplain}

\end{document}